\documentclass[12pt]{amsart}

\newcommand {\supplus}{\mathop{{\supset}\llap{\raise
0.5pt\hbox{\normalfont\small+}\hskip 0.5pt}}}

\newcommand {\subplus}{\mathop{{\subset}\llap{\raise
0.5pt\hbox{\normalfont\small+}\hskip 0.5pt}}}

\newcommand {\Cee}    {{\mathbb  C}}

\newcommand {\Zee}    {{\mathbb  Z}}

\newcommand {\fder}   {{\mathfrak{der}}}   %

\newcommand {\fg}     {{\mathfrak{g}}}    %
\newcommand {\fgl}    {{\mathfrak{gl}}}  %
\newcommand {\fh}     {{\mathfrak{h}}}

\newcommand {\fk}     {{\mathfrak{k}}}

\newcommand {\fN}     {{\mathfrak{N}}}  %
\newcommand {\fn}     {{\mathfrak{n}}}

\newcommand {\fpo}    {{\mathfrak{po}}}

\newcommand {\fsl}    {{\mathfrak{sl}}}

\newcommand {\fsp}    {{\mathfrak{sp}}}

\newcommand {\fsvect} {{\mathfrak{svect}}}

\newcommand {\fvect}  {{\mathfrak{vect}}}   %

\newcommand {\fwitt}  {{\mathfrak{witt}}}

{\newcommand {\cal} {\mathcal}

\def \opname#1#2%
  {\expandafter\newcommand \csname #1\endcsname {{\mathop{#2}\nolimits}}}

\def\rmname#1%
  {\expandafter\newcommand \csname #1\endcsname {{\operatorname{#1}}}}

\rmname{act}
\rmname{Ad}
\opname{ad}     {{ad}}
\rmname{Alt}
\rmname{alt}
\rmname{Ann}
\opname{antidiag} {{antidiag}}
\opname{Ber}    {{Ber}}
\rmname{ber}
\rmname{card}
\rmname{Char}
\rmname{cem}
\rmname{cj}
\rmname{cntr}
\rmname{codim}
\rmname{const}
\rmname{col}
\rmname{cpr}
\opname{diag}{{diag}}
\opname{Der}{{Der}}
\opname{End} {{End}}
\rmname{gr}
\rmname{Hom}
\rmname{Ht}
\rmname{Id}
\rmname{id}
\rmname{ind}
\rmname{Le}
\rmname{mult}
\opname{Mor} {{Mor}}
\rmname{nm}
\opname{Ob} {{Ob}}
\rmname{Pic}
\rmname{pr}
\rmname{pro}
\rmname{Prime}
\rmname{Proj}
\rmname{prt}
\rmname{Q}
\rmname{qet}
\rmname{qtr}
\rmname{rd}
\rmname{rk}
\rmname{row}
\rmname{Res}
\rmname{scalar}
\opname{Ser} {{Ser}}
\rmname{Smbl}
\rmname{str}
\rmname{sgn}
\rmname{sq}
\rmname{symm}
\rmname{supp}
\rmname{St}
\opname{Spec} {{Spec}}
\rmname{Spm}
\rmname{tr}
\rmname{vpt}

\opname{vol}  {{v\hspace{-0.1ex}o\hspace{-0.02ex}l\/}}
\opname{pnt}  {\text{\normalfont pt}}
\opname{Span} {{Span}}
\opname{slim} {\overline{\lim}}
\opname{Vol}  {{V\hspace{-0.55ex}o\hspace{-0.02ex}l\/}}
\opname{Par}  {{P\hspace{-0.3ex}a\hspace{-0.05ex}r\/}}

\rmname{Mat}
\rmname{Bil}
\rmname{Diff}
\rmname{Ker}
\rmname{Herm}
\rmname{Coker}
\rmname{Conn}
\rmname{Covect}
\rmname{Vect}
\rmname{Int}

\opname{IM} {\text{\normalfont Im}}
\opname{RE} {\text{\normalfont Re}}

\opname{Aut} {{A\hspace{-0.2ex}u\hspace{-0.1ex}t\/}}
\opname{GL} {{G\hspace{-0.3ex}L}}
\opname{SL} {{S\hspace{-0.3ex}L}}
\opname{Exp} {{E\hspace{-0.2ex}x\hspace{-0.1ex}p\/}}
\opname{GQ} {{G\hspace{-0.2ex}Q}}
\opname{OSp} {{O\hspace{-0.25ex}S\hspace{-0.15ex}p\/}}
\opname{Out} {{O\hspace{-0.25ex}u\hspace{-0.15ex}t\/}}
\opname{Spp} {{S\hspace{-0.2ex}p\/}}
\opname{SpO} {{S\hspace{-0.2ex}p\hspace{-0.02ex}O\/}}
\opname{Pe} {{P\hspace{-0.25ex}e\/}}
\opname{SPe} {{S\hspace{-0.25ex}P\hspace{-0.25ex}e\/}}
\opname{Spin} {{S\hspace{-0.25ex}p\hspace{-0.05ex}i\hspace{-0.1ex}n\/}}
\opname{Iso} {{I\hspace{-0.25ex}s\hspace{-0.1ex}o\/}}
\opname{SSPe} {{S\hspace{-0.25ex}S\hspace{-0.15ex}P\hspace{-0.25ex}e\/}}
\opname{PeU} {{P\hspace{-0.25ex}e\hspace{-0.1ex}U\/}}
\opname{QU} {{Q\hspace{-0.15ex}U\/}}
\opname{U} {{U\/}}

\opname{cGQ} {{\cal G \hspace{-0.2em} Q \/}}
\opname{cSL} {{\cal S \hspace{-0.2em} L \/}}
\opname{cGL} {{\cal G \hspace{-0.2em} L \/}}
\opname{cOSp} {{\cal O \hspace{-0.2em} S \hspace{-0.3em} \it p\/}}
\opname{cPe} {{\cal P \hspace{-1.5pt} \it e\/}}
\opname{cVect} {{\cal V \hspace{-1.5pt} \it
e\hspace{-0.1ex}c\hspace{-0.1ex}t\/}}
\opname{cVol} {{\cal V \hspace{-1.5pt} \it o\hspace{-0.1ex}l\/}}
\opname{cAut} {{\cal A \hspace{-0.2em} \it u\hspace{-0.1em}t\/}}
\opname{cCovect} {{\cal C \hspace{-1.5pt}
     \it o\hspace{-0.1ex}v\hspace{-0.1ex}e\hspace{-0.1ex}c\hspace{-0.1ex}t\/}}
\opname{CW} {{C\hspace{-0.15ex}W}}

\newcommand {\eps} {\varepsilon}

\newcommand {\pder}[1] {{\frac{\partial}{\partial {#1}}}}

\newbox\tmpbox
\def\vect#1%
{ \setbox\tmpbox = \hbox{$#1$}
  \ifdim \wd\tmpbox<2ex
     \vec{#1}
  \else
     #1\llap{\overrightarrow{\phantom{#1\!\!}\,}}
  \fi
}

\newcommand {\bcdot}   {\mathbin{\hbox{\raise.4ex\hbox{\bf.}}}} 

\def \nopoint#1#2{}

\newcommand {\secno} {}

\newtheorem{Theorem}{\secno Theorem}

\newenvironment {th*}[1]
    {\gdef\thname{#1} \begin{thn}}%
    {\end{thn}}
\newtheorem{thn}[Theorem] {\thname}

\theoremstyle{definition}

\newenvironment {ex*}[1]
    {\gdef\thname{#1} \begin{exn}}%
    {\end{exn}}
\newtheorem{exn}[Theorem]{\thname}

\theoremstyle{remark}

\newenvironment {rem*}[1]
    {\gdef\thname{#1} \begin{remn}}%
    {\end{remn}}
\newtheorem{remn}[Theorem]{\thname}

\newcommand {\ssec}{\subsection*}

\setcounter{tocdepth}{1}

\begin{document}

\title[DEFINING RELATIONS FOR CLASSICAL LIE ALGEBRAS]{DEFINING RELATIONS FOR
CLASSICAL LIE ALGEBRAS\\ OF POLYNOMIAL VECTOR FIELDS}

\author{Dimitry Leites, Elena Poletaeva}

\address{Department of Mathematics,
Stockholm University, Roslagsv.  101, S-104 05, Stockholm, Sweden
(for correspondence) mleites@math.su.se; Department of
Mathematics, University of California at Riverside, Riverside,
USA; elena@math.ucr.edu}



\begin{abstract}
We explicitly describe the defining relations for simple Lie algebra
$\fvect (n)= \fder~\Cee [x]$, where $x=(x_1, \dots, x_n)$, of vector
fields with polynomial coefficients and its subalgebras of divergence
free, hamiltonian and contact vector fields, and for the Poisson
algebra (realized on polynomials).  We consider generators of $\fvect
(n)$ and its subalgebras corresponding to the system of simple roots
associated with the standard grading of these algebras.  (These
systems of simple roots are {\it distinguished} in the sense of
Penkov--Serganova \cite{PS}.)
\end{abstract}

\maketitle

\section*{Introduction}
This is our paper published in Mathematica Scandinavica, 81, 1997,
5--19. We just wish to make it more accessible. References are
updated.

The class of simple $\Zee$-graded Lie algebras of polynomial
growth (SZGLAPGs, for short) over $\Cee$ often appears in various
problems of mathematics and physics.  V.~Kac conjectured \cite{K}
(and O.~Mathieu proved this \cite{M}) that this class consists of
the following algebras:

1) finite dimensional Lie algebras (each determined by a Cartan matrix
or Dynkin diagram);

2) (twisted) loop algebras (each determined by the Cartan matrix
corresponding to an extended Dynkin diagram);

3) four series of Lie algebras of vector fields with polynomial
coefficients or, briefly, \lq\lq vectorial" algebras.
(\footnotesize The name \lq\lq vectorial" is influenced by
physical terminology. Some physicists working with string theories
call the Lie superalgebras pertaining to these theories \lq\lq
stringy".  This term is very accurate and doubly suggestive: it
reminds the relation of these algebras to string theories and
reveals that the stringy algebras do look as a bunch of strings
--- modules over the Witt algebra.  Stringy superalgebras are a
particular case of what we propose to call \lq\lq vectorial
algebras": with Laurent polynomials as coefficients.  We concede
that \lq\lq vectorial" might sound a bit reckless; still, in our
terminology there is no chance to look at something like a \lq\lq
Cartan subalgebra of a Lie subalgebra of Cartan type with Cartan
matrix".\normalsize)

$\fvect (n) = \fder~\Cee [x]$, where $x=(x_1, \dots, x_n)$, the {\it
general vectorial algebra};

$\fsvect (n) = \{D \in \fvect (n): \text{div }D = 0\}$,   the
{\it special or divergence-free algebra};

$ \fh (2n)  = \{D \in {\fvect}(2n):\ D\omega = 0\}$ for $\omega  =
\sum_i dp_i\wedge dq_i$, here $x=(q, p)$,   the {\it
Hamiltonian  algebra};

$\fk (2n+1)  = \{D \in  \fvect (2n+1):\ D\alpha = f\alpha \}$\
for $ \alpha =  dt +\sum_i (p_idq_i + q_idp_i)$, here $x=(t, q, p)$,  the
{\it contact} Lie algebra;

4) this class consists of one
\lq\lq stringy" (i.e., pertaining to the string theory) Lie algebra: the {\it
Witt algebra}
$\fwitt = \fder \Cee[x^{-1}, x]$.

\ssec{A central extension} The Lie algebra $\fpo (2n)$ whose space is
$\Cee [q, p]$ and the bracket is the Poisson bracket is called the
{\it Poisson algebra}.  It is the nontrivial central extension of $\fh
(2n)$, with the 1-dimensional center generated by constants.  Its
geometric interpretation: $\fpo (2n)$ preserves the connection with
the above form $\alpha$ in a line bundle over $\Cee^{2n}$ with a
symplectic structure.  Though $\fpo (2n)$ is not simple, it is useful
in applications (similarly, Kac--Moody algebras are more useful than
the loop algebras).  So we will consider it as well.

The end product of the deformation (the physicists call it \lq\lq {\it
quantization}") of $\fpo (2n)$ is \lq\lq $\fgl" (\Cee [q])$.  (We used
quotation marks because for infinite dimensional spaces there are many
distinct \lq\lq $\fgl$"s, cf.  \cite{E}.)  Our results might help to
understand how to quantize and \lq\lq $q$-quantize" the Poisson
algebra of functions on $2n$-dimensional torus.  In particular, our
results are at variance with expectations of \cite{S}.

\ssec{Related results} The description of the above algebras in terms
of defining relations is vital in questions where it is necessary to
identify an algebra determined by its defining relations
(Eastbrook-Wahlquist prolongations, Drinfeld's quantum algebras,
etc.).  Presentations of the algebras from the first two classes are
known.  We will give a presentation for the remaining of the above
algebras.

It was in \cite{FF1} that Feigin and Fuchs first published a presentation
of a nilpotent part ($\fg_+$ in our notations, see sect.  0.3) of
$\fg=\fvect (n)$, later generalized (with a gap in the case of
$\fg=\fh(2n)$) for all simple vectorial Lie algebras in \cite{FF2}.
However, both the generators and relations considered in \cite{FF1} and
\cite{FF2} are too numerous.  Besides, they were implicit.

This was no wonder: in a similar situation for loop and affine
Kac--Moody algebras $\fg^{(1)}=\fg\otimes\Cee[t^{-1}, t]$ with the
grading $\deg g=0$ for $g\in\fg$, $\deg t=1$ it is more advisable to
consider the generators of $\fg$ and the lowest and highest generators
of $\fg\otimes t$ and $\fg\otimes t^{-1}$, respectively.  So we
suggest to enlarge the algebra considered by Feigin and Fuchs.  In the
same manner as is done for loop or Kac--Moody algebras, we add to
$\fg_+$ the positive root vectors of $\fg_0$.

In this way we only have to slightly work on the result of \cite{FF1}
and \cite{FF2} to get an explicit answer.  Our results are based on
the general theorems from \cite{FF1}, \cite{FF2} related to the
Hochschield--Serre's spectral sequence applied to locally nilpotent
Lie algebras.  (Similar statements for Kac--Moody algebras are
well-known, cf.  \cite{K}.)  The exceptional cases of small number of
indeterminates, not covered by the general theorems of \cite{FF1},
\cite{FF2} were studied in \cite{HP1}, \cite{HP2} by means of a
computer.

\footnotesize

This paper, preprinted in \cite{L}, no.  31, is the first in a series
of papers devoted to presentation of simple $\Zee$-graded Lie {\it
super}algebras of polynomial growth over $\Cee$: the main ideas and
some examples were first delivered by D.~Leites at L.~Faddeev's
seminar in 1981, and later at the lectures at CWI, Amsterdam in 1986.

The publication of the text was delayed for 10--15 years for no
reason; meanwhile there appeared several results with description of
presentations of simple Lie algebras and Lie superalgebras, or their
bases, connected with this paper: \cite{GL1}, \cite{HP1}, \cite{HP2},
\cite{LSe}.

A.~Dzhumadildaev informed us that, with a student, he recently
(1992) obtained a description \cite{DI} identical to ours for
simple vectorial Lie algebras over an algebraically closed field
of large characteristics ($p>7$) and for the simplest value of the
parameter $\underline{N}$, namely, for $(1, \dots, 1)$; (for
definitions, see \cite{GL}).

\normalsize

\section*{\protect\S 0. Preliminaries}

\ssec {0.1.  Defining relations} Presentations of Lie algebras of
classes 1) and 2) of the above list of SZGLAPGs are known; the most
popular is the one that can be neatly encoded in terms of a graph
(Dynkin diagram).  Still, we have to explain what does this common
knowledge amount to since there is no \lq\lq natural" set of
generators for a simple Lie algebra.

The common presentation of the simple Lie algebras is performed in
terms of Chevalley generators.  A rival, less popular but also
very useful in various problems, set of generators for any finite
dimensional simple Lie algebra is to pick just two generators (as
suggested first by N.~Jacobson).  These generators are related to
the principal embedding of $\fsl(2)$; for a proof, see \cite{BO},
and for applications \cite{F2}, \cite{FNZ}, \cite{LS}.  What
should be meant under {\it relations} for such a choice of
generators and the degree of the ambiguity of these relations is,
however, very far from clear, cf.  \cite{F2}, \cite{FNZ}; for an
answer, see \cite{GL1}.

Contrarywise, for any nilpotent Lie algebra $\fn$, the problem has
a natural and unambiguous solution: a basis of the space $\fn
/[\fn, \fn ] = H_{1}(\fn )$ is a set of {\it generators} of $\fn
$. Suppose, as will be the case in our examples, that there is a
set of outer derivations acting on $\fn$ so that $\fn /[\fn , \fn
]$ splits into the direct sum of 1-dimensional eigenspaces.  Then
the choice of a basis is unique up to scalar factors.

To describe relations between the generators, consider the standard
homology complex for $\fn $ with trivial coefficients (\cite{Fu}):
$$
0 \longleftarrow \fn \stackrel{d_{1}}{\longleftarrow} \fn \wedge \fn
\stackrel{d_{2}}{\longleftarrow} \fn \wedge \fn \wedge \fn
\stackrel{d_{3}}{\longleftarrow} \dots \, .
$$
By definition
$$
d_{1}(x\wedge y)=[x, y],\quad d_{2}(x\wedge y\wedge z)=[x, y]\wedge z
+ [y, z]\wedge x + [z, x]\wedge y.
$$
The elements of $\Ker d_{1}$ are, obviously, relations.  The
consequences of the Jacobi identity should be considered as
trivial relations; they constitute $\IM\; d_{2}$.  Thus, a basis
of $H_{1}(\fn)$ allows to construct the generators of
$H_{2}(\fn)$, and applying $d_1$ to them we obtain {\it defining
relations} for a given nilpotent Lie algebra $\fn$.  Observe that
the same arguments apply as well to {\it locally nilpotent} Lie
algebras (such as the algebras $\fN _{\pm}$ and $\fg_{\pm}$
considered in \cite{FF1}, \cite{FF2} and below).

\ssec {0.2.  Serre's relations} It was J.-P. Serre who, perhaps, for
the first time, wrote down relations between Chevalley generators of a
simple finite dimensional Lie algebra $\fg$.  They can be also written
for the algebras of class 2) and more general Lie algebras with Cartan
matrix, cf.  \cite{GKLP} and \cite{K}; \cite{GL2}.

The relations are expressed in terms of Chevalley generators, so let
us first define Chevalley generators.  Let
$$
{\fg} ={\fn_{-} }
\oplus \fh \oplus \fn _{+}
$$
be the root decomposition of $\fg$, where $\fh$ is a {\it maximal
torus} or, more exactly, toral, i.e., diagonalizable, subalgebra and
$\fn_{\pm}$ are maximal nilpotent subalgebras generated by the root
vectors corresponding to all positive (for $+$), respectively negative
(for $-$), roots.  (\footnotesize The notion of a Cartan subalgebra
--- a nilpotent subalgebra coinciding with its normalizer --- plays a
modest role in the theory of Lie algebras over fields of prime
characteristics, or for superalgebras, or in the case we are
considering, namely, of infinite dimension, or combination of
these classes.  A reason: for the simple finite dimensional Lie
algebras over $\Cee$, the Cartan subalgebra is both (maximal)
commutative and diagonalizable.  Sometimes, e.g., for
$\fpo{(2n)}$, $\fh{(2n)}$, it is neither.\normalsize ) By setting
$$
\deg ( X_{\pm \alpha})= \pm 1
$$
for a root vector corresponding to a simple root $\alpha$ (or its
opposite) we endow $\fn_{\pm}$ with a $\Zee$-grading, and it turns out
that
$$
[\fn_{\pm},\fn_{\pm}]= \{X\in{\fn_{\pm}};\; \; \deg X \geq 2
\quad (\text{for}\quad +),\quad \deg X \leq -2 \quad (\text{for}\quad -) \}.
$$
Hence, the elements $X_{\pm \alpha}$ themselves represent their
homology classes in $H_{1}(\fn_{\pm})$, namely, $H_{1}(\fn_{\pm}) =
\mathop{\oplus}\limits_{\alpha_i \in S} \Cee X_{\pm \alpha _i},$ where
$S$ is any system of simple roots.  Informally speaking, these
generators are {\it pure}: they are not defined modulo anything.

Let us normalize these generators as follows.  Set $X_{i}^{\pm}=X_{\pm
\alpha_{i}}$ and introduce auxiliary generators $H_i = [X^{+}_i,\
X^{-}_i]$.  Using rescaling $H_i \mapsto \lambda H_i$ it is possible
to select the $X^{\pm}_i$ so that
$$
[H_i,\ X^{\pm}_j] = \pm A_{ij}X^{\pm}_j,
$$
where the matrix $(A_{ij})$ is normed so that $A_{ii} = 2$ for all
$i$.  The normed matrix $(A_{ij})$ is called the {\it Cartan matrix}
of $\fg$, the generators $X^{\pm}_i$ of $\fn^{\pm}$ the {\it Chevalley
generators}.

The {\it relations} between Chevalley generators $ X_{i}^{\pm}$ of
$\fn_{\pm}$ are called {\it Serre's relations}.  They are expressed in
terms of the Cartan matrix as follows:
$$
(ad(X_{i}^{\pm}))^{1-{A_{ij}}} (X_{j}^{\pm})=0 \quad (i\not= j).
\eqno{(SR_{\pm})}
$$
Serre's relations represent homology classes from $H_{2}(\fn_{\pm})$
which are also pure:
$$
H_{2}(\fn_{\pm}) = \mathop{\oplus}\limits_{\alpha _i \not=
\alpha _j \in S}\Cee
(X_{\pm \alpha _i}\wedge X_{\pm r_{\alpha _i}(\alpha _j)}),
$$
where $r_{\alpha _i}$ is the reflection associated to the root $\alpha
_i$: $r_{\alpha _i}(\alpha _j) = \alpha _j - A_{ij}\alpha _i$.  The
elements $X_{\pm \alpha _i}\wedge X_{\pm r_{\alpha _i}(\alpha _j)}$
correspond (via $d_1$) to the relations $(SR_{\pm})$.

\begin{rem*}{Remark}
The Serre's relations show that the lowest weight of the
$\fsl(2)$-module generated by $X^{+}_j$ (or the highest weight of the
module generated by $X^{-}_j$) is equal to $A_{ij}$ (resp.  $-A_{ji}$)
with the copy of $\fsl(2)$ we are talking about being generated by
$X^{\pm}_i$.
\end{rem*}

There are also relations between $\fn_{+}$ and $\fn_{-}$ (and\ $\fh$):
$$
[H_{i},H_{j}]=0,\;\;
[X_{i}^{+},X_{j}^{-}]= \delta_{ij}H_{i},\;\;
[H_{i},X_{j}^{\pm}]=\pm{A_{ij}}X_{j}^{\pm}. \eqno{(SR_{0})}
$$

One of the reasons why the above relations (SR) won priority over any
other presentation is that they can be encoded in a single
nice-looking graph, the {\it Dynkin diagram} of ${\fn}^{\pm}$ (or,
which is the same, of $\fg$), see \cite{B}, \cite{OV}.  There are, however,
other possibilities, and recently it became clear that they are not so
awful, cf.  \cite{GL1}, \cite{GL2}.

\ssec{ 0.3.  Vectorial algebras} About 1980 D.~Leites conjectured
that the Lie algebras of vector fields with polynomial
coefficients are finitely presented.  I.~Kantor informed us then
that much earlier he also arrived to the same conjecture; he
believed that these algebras are determined by a pair of \lq\lq
Cartan matrices" or \lq\lq Dynkin diagrams" and even produced
diagrams hypothetically corresponding to ${\fvect}(n)$.  As we
will see Kantor's conjecture was almost true but definitely not
that simple, as is clear from \cite{GL2}, especially the last
section.

A weaker conjecture, on finite determinacy, had been independently
proved by Bondal and Ufnarovsky in \cite{BU} as a by-product in their
answer to another problem.  They did not appreciate, however, the
importance of this by-product, and no explicit formulas were written
neither then nor later, cf.  \cite{U}.

By the time \cite{BU} was out of print it became manifest that the
simplest examples shatter any hope for neat relations.  Indeed,
consider the $\fn_{\pm}$-part of the very first examples, $\fwitt$ and
${\fvect}(1)$.  Recall that a natural basis in $\fwitt$ is
$$
e_{i}=x^{i+1}\frac{d}{dx}\quad \text{ for}\; i\in\Zee
$$
with relations
$$
[e_{i},e_{j}]=(j-i)e_{i+j}.
$$
Clearly, $\Span(e_{0})$ is the maximal torus, $\fn_{\pm}$ is generated
by $e_{\pm1}$ and $e_{\pm2}$.  (Note that $\fn_-$ for $\fvect (1)$ is
just $\Span(e_{-1})$.)

The relations in $\fn_{-}$ and $\fn_{+}$ and between them, not as
popular as Serre's relations, are, nevertheless, known, thanks to
B.~Feigin \cite{F1}, since 1979, and later were rediscovered
(\cite{FNZ}):
\[
\begin{aligned}
\begin{matrix}
[e_{-1}, e_{1}]=2e_{0},\; [e_{-2}, e_{2}]=4e_{0},\; [e_{0}, e_{i}]=ie_{i}
\; (i=\pm 1, \pm 2),\\{}
[e_{-1}, e_{2}]=3e_{1},\; [e_{1}, e_{-2}]=-3e_{-1};
\end{matrix}
\qquad (DR_{0}) \\
(\ad e_{\pm 1})^{3}e_{\pm 2} + 6(\ad e_{\pm 2})^{2} e_{\pm 1}=0,\quad
(\ad e_{\pm 1})^{5}e_{\pm 2} +40 (\ad  e_{\pm 2})^{3}e_{\pm 1}=0.
\qquad (DR_{\pm})
\end{aligned}
\]
The relations that do not involve $e_{-2}$ are the defining relations for
$\fvect(1)$.

The relations of degree ${\pm 7}$ are \lq\lq dirty".  This means that
(say, for the plus sign) in the 2-dimensional space of cycles
$$
\Span(3e_1\wedge e_6 -5e_2\wedge e_5,\; \;  e_2\wedge e_5
-3e_3\wedge e_4)
$$
the relations span any line transversal to the line of boundaries:
$$
d_{3}(e_1\wedge e_2 \wedge e_4) = e_3\wedge e_4 -3e_5\wedge e_2
+3e_6\wedge e_1.
$$

By that time (1980) D.~Fuchs himself got interested in the problem
and made a major step forward.  He studied the relations for
$\fg_+=\mathop{\oplus}\limits_{i>0}\fg _i$, where $\fg$ is a Lie
algebra of vector fields (class 3) above) in the standard
$\Zee$-grading ($\deg x_i = 1$ for all $i$, except for $\deg t =
2$). For $\fsvect (2) = \fh (2)$ and for $\fk(3)$, there are also
some dirty relations.

Fuchs showed that, for vectorial Lie algebras depending on
sufficiently large number of indeterminates, there are no
``dirty'' relations and this \lq\lq large" number is actually
pretty small: it is equal to 3 except for the series $\fk$ when it
is equal to 5.

Probably overjoyed with the above result, Fuchs conjectured
(\cite{FF2} and references therein) that for any simple vectorial
$\fg$ in the standard grading all the relations between the generators
of $\fg_+$ are of degree 2, namely,
$$
H_{2}(\fg_{+})=[H_{2}(\fg_{+})]_ {2} =
\wedge ^{2}(\fg_{1})/\fg_{2}.
$$
Fuchs proved this conjecture for $\fg = \fvect (n)$. The statement
and its proof (\cite{FF1}) is also true for the other series of
simple infinite dimensional vectorial Lie algebras except for the
Hamiltonian algebras for which we provide with the correct
statement.  Fuchs' statement might have been corrected already in
1984 when Kochetkov discovered that the similar relations for the
Lie superalgebra $\fh (0|n)$ always contain a component of degree
3; we should have put more faith in Kochetkov's result from the
very beginning, cf. \cite{GKLP}, \cite{HP1}, \cite{HP2}.

These relations, even correct ones, are, however, implicit and
inconvenient: if $\dim \fg_{-1} = n$, then there are
$\dim\fg_1\approx n^3$ generators of $\fg_+$ and $\approx n^6$
relations between these generators.

In what follows we will proceed precisely as for the (twisted)
loop or Kac-Moody algebras, namely, replace $\fg_+$ with the
maximal locally nilpotent subalgebra $\fN_{+}$ of $\fg$: this
steeply diminishes the number of generators and relations,
moreover, the relations become graphic.

Thus, we look at the 2nd term of Hochschield-Serre's spectral
sequence constructed for the pair $\fg_{\pm} \subset \fN_{\pm}$
(see \S 2) and unite the bases of three spaces obtained
($H_i(\fn_{\pm}$; $H_j(\fg_{\pm}))$ for $i+j=2$, where $\fg_{+} =
\mathop{\oplus}\limits_{i>0}\fg_i$ and $\fg_{-}=
\mathop{\oplus}\limits_{i<0}\fg_i$); the first two of them are
calculated with the help of the Borel--Weil--Bott theorem, to the
third one we apply (amended) Fuchs' result. These three spaces
give an upper bound for the space of relations. More accurate
study (of the behavior of the codifferential) may diminish the
number of independent relations; for such accurate analysis, see
\cite{P}.

We started this work with expectations to get simple, Serre-type
looking relations.  To an extent our relations {\it are} simple if we

a)  do NOT express them in terms of generators only (likewise,
Chevalley generators contain \lq\lq redundant" generators $H_i$
expressible in terms of $X_i^{\pm}$);

b) simplify them by taking not necessarily the one which corresponds
to the lowest weight of the corresponding $\fg_0$-module but a linear
combination of it and other relations of the same weight obtained from
other modules of relations in order to get a simpler, \lq\lq more
factorizable", cycle.

\ssec{0.4.  The main result} We give an explicit expression for the
generic case of defining relations in the maximal locally nilpotent
subalgebras $\fN_{\pm} = \fn_{\pm}\oplus\fg_{\pm}$ generated by
positive and negative root vectors of $\fg$ in the standard
$\Zee$-grading of $\fg$, where $\fg_{\pm}$ are described above and
$\fn_{\pm}$ are the maximal nilpotent subalgebras of $\fg_0$ described
in textbooks (e.g., \cite{B}, \cite{OV}).

The exceptional cases of small dimension are investigated by
N.~Hijligenberg and G.~Post, for explicit formulas, see
\cite{HP1}, \cite{HP2}.

\ssec{0.5.  A problem: what is an analog of the Weyl group for a
given vectorial algebra} Unlike the case of simple finite
dimensional Lie algebras with Cartan matrix over $\Cee$, there are
several nonisomorphic types of maximal nilpotent subalgebras
$\fg_+$ of a vectorial Lie algebra.  (The choice we made is
distinguished by the fact that the subalgebra of elements of
nonnegative degree is the unique maximal subalgebra of finite
codimension; it also happens that it is easier to deal with.)

We encounter similar phenomenon with (say, simple finite
dimensional) Lie superalgebras, even if they have Cartan matrices,
or in prime characteristic.  The experience with Lie superalgebras
(\cite{LSe}, \cite{FLV}) suggests to consider all maximal
nilpotent subalgebras (or, equivalently, bases, i.e., systems of
simple roots) simultaneously and consider the group that
transitively acts on the bases --- an analog of the Weyl group.  A
solution to this problem will be discussed elsewhere, see recent
\cite{SV}.

\section*{\protect\S 1. Generators in vectorial Lie algebras}
Set $\partial _i= \pder{x_i}$.

\ssec{1.1.  Generators of $\fvect (n)$} Some of the generators of
$\fvect (n)$ generate its following subalgebras:

\small
\begin{tabular}{|c|c|c|}
\hline
&$\fsl (n+1)$&\\
\hline
$\fN_+$& $x_1\partial _2,  \; \; \; \; \; \;   \dots ,  \; \; \;
x_{n-1}\partial  _n, \; \; \;   \pmb{x_n\sum x_i\partial _i}$ &
$x_n^2\partial _1$\\
\hline
$\fN_-$&$x_2\partial _1, \; \; \;  \dots , \; \; \;
x_n\partial  _{n-1}, \; \;
\; \; \partial  _n$&\\
\hline
\text{ notations}& $X^{\pm}_1,\; \; \; \;    \dots , \; \; \;
X^{\pm}_{n-1},
\; \; \; X^{\pm}_n$&$Y$\\
\hline
\end{tabular}

\vskip 0.5cm

\normalsize

The generators of $\fsvect (n)$ are the same as of $\fvect (n)$ but
without the boldfaced element $X^+_n$.

\ssec{1.2.  Generators of \protect $\fk (2n+1)$} First of all it is
convenient to express the elements of $\fk (2n+1)$ and its subalgebras
$\fpo (2n)$ and $\fh (2n)$ in terms of {\it generating functions} (for
the series $\fh (2n)$ these generating functions are called {\it
Hamiltonians}) as follows.

To every $f\in \Cee[q, p]$ assign the Hamiltonian field $H_f$
corresponding to the Hamiltinian $f$:
$$
H_f=\sum\biggl(\frac{\partial  f}{\partial  p_i}\frac{\partial  }
{\partial  q_i}\ -\ \frac{\partial  f}{\partial  q_i}\frac{\partial  }
{\partial  p_i}\biggr).
$$
To every $f\in \Cee[t, q, p]$ assign the contact field
$$
K_f=\triangle (f)\pder{t}+ \frac{\partial f} {\partial  t}E+
H_f,
$$
where $\triangle (f)= 2f - E(f)$ for $E=\sum_i y_i \frac{\partial
}{\partial y_i}$; here the $y$ are all the coordinates except $t$.

In particular, the functions that not depend on $t$ generate the Poisson
algebra realized by vector fields $K_f$.

To the commutator of hamiltonian vector fields $[H_f, H_g]$ there
corresponds the {\it Poisson bracket} of generating functions, and to
the commutator of contact vector fields $[K_f, K_g]$ there corresponds
the {\it Lagrange} or, as it is more often called, {\it contact}
bracket:
\begin{eqnarray}
\{f,g\}_{c.b.}=\triangle (f)\frac{\partial  g}{\partial  t}\ -\ \frac{\partial
f} {\partial  t}\triangle (g)-\{f,g\}_{P.b.},\nonumber
\end{eqnarray}
and the Poisson bracket is given by the
formula
\begin{eqnarray}
\{f,g\}_{P.b.}=\sum\biggl(\frac{\partial  f}{\partial  p_i}\frac{\partial  g}
{\partial  q_i}\ -\ \frac{\partial  f}{\partial  q_i}\frac{\partial  g}
{\partial  p_i}\biggr)\nonumber
\end{eqnarray}

It is well-known that
$$
\fh (2n) = \Span(H_f: f\in \Cee[q, p]),\quad  \fk (2n+1) = \Span(K_f:
f\in \Cee[t, q, p]).
$$

In what follows we will by abuse of language write just $f$ instead of
$H_f$ or $K_f$; in so doing we must remember that the degree of the
vector field $K_f$ or $H_f$ generated by a monomial $f$ of degree $k$
is equal to $k-2$.

The important for us subalgebras of $\fk (2n+1)$ and their generators
with our notations are:

\vskip 0.5 cm

\small
\begin{tabular}{|c|c|c|}
\hline
&$\fsp (2n+2)$&\\
\hline
$\fN _- $& $p_1\; \; \; \;      q_1p_2\; \; \;    \dots \; \; \;
q_{n-1}p_n\; \;
\; \;   q_n^2$&\\
\hline
$\fN _+$ & $\pmb{tq_1}\; \; \; \;   p_1q_2 \; \;  \dots \; \;  q_np_{n-1} \;
\;
\; \; p_n^2 $&$q_1^3$\\
\hline
\text{ notations} &$X^{\pm}_0,\; \; \;  X^{\pm}_1,\; \;  \dots \; \;
X^{\pm}_{n-1},\; \; \;  X^{\pm}_n $&$Y$\\
\hline
\end{tabular}
\normalsize
\vskip 0.5 cm

The generators of $\fh (2n)$ and $\fpo(2n)$ are those above without
the boldfaced element $X^+_0 = tq_1$.

\section*{\protect\S 2. How to calculate relations}
Clearly, $\fN_-$ for $\fvect(n)$ and $\fsvect (n)$ coincides with
$\fn_-$ for $\fsl(n+1)$; $\fN_-$ for $\fk (2n+1)$ and $\fpo (2n)$
coincides with $\fn_-$ for $\fsp(2n+2)$.  Therefore, the relations for
these Lie algebras $\fN_-$ are known.

The remaining defining relations are found with the help of the
Hochschield-Serre spectral sequence \cite{Fu} for the pair
$\fg_{\pm} \subset \fN_{\pm}$ and the results of \cite{FF2} on
$H_{2}(\fg_{+})$.  We estimate the relations from above: the space
of relations is contained in the direct sum of the following
spaces:
$$
\begin{matrix}
H_{2}(\fn_{\pm});\hfill (1)\\
H_{1}(\fn_{\pm}; H_{1}(\fg_{\pm }))\hfill (2)\\
\begin{matrix}H_{0}(\fn_{+}; H_{2}(\fg_{+}))=\Span(\text{the } \fn
_{-}\text{-lowest vectors of}\ H_{2}(\fg_{+}));\hfill (3^+)\\
H_{0}(\fn_{-}; H_{2}(\fg_{-}))  =\Span(\text{the } \fn _{+}\text{-highest
vectors of }\ H_{2}(\fg_{-})).\hfill (3^-)
\end{matrix}
\end{matrix}
$$
Homology (1) are known, cf.  \cite{OV}; homology (2) are calculable
with the help of the Borel--Weil--Bott theorem, cf.  \cite{Go}.  When
homologies (3) are known, it remains to compute the differential in
the spectral sequence.  Usually it is zero.

Notice that the first two types of these relations
are of the form
$$
(\ad  X_i^{\pm})^{k_{ij}}(Y_j^{\pm}) = 0,
$$
where the $X_i^+$ are the generators of $\fn _+$ and the $ Y_j^+$ are
the lowest weight vectors of the $ \fg _0$-module $ \fg_1$; the $
X_i^-$ are the generators of $\fn _-$ and the $ Y_j^-$ are the highest
weight vectors of the $\fg _0$-module $ \fg_{-1}$.

\begin{rem*}{Remark} Similar calculations by induction on the rank for
simple finite dimensional and (twisted) loop algebras give the
shortest known to us proof of completeness of the Serre defining
relations.
\end{rem*}

\section*{\protect\S 3.  Relations in $\fk(2n+1)$, $\fpo(2n)$ and
$\fh(2n)$}

\ssec{3.1.  Relations in \protect $\fN_-$ for $\fk(2n+1)$ and
$\fpo(2n)$ } These relations are the same as for $\fn_-$ of
$\fsp(2n+2)$.

\ssec{3.2.  Relations in $\fN_-$ for $\fh(2n)$} The Lie algebra
$\fh(2n)$ is generated by the same elements as $\fk(2n+1)$ and
$\fpo(2n)$ but the relations are different: for $\fh(2n)$ there is an
additional relation of weight $(0, \dots, 0)$ with respect to
$\fsp(2n)$ because
\[
H_2(\fg_{-1}) = \begin{cases}\Lambda^2(\fg_{-1}) = R(\pi_2) \oplus
R(0)&\text{if \(n>1\)}\\
R(0)&\text{if \(n=1\)},\end{cases}
\]
where $R(\pi_2)$ and $R(0)$ denote the irreducible $\fsp(2n)$-modules
with highest weights $\pi_2$ and $0$, respectively; $\pi_i$ is the
$i-th$ fundamental weight (see \cite{OV}).  For $n = 1$ the
corresponding relation is of the form
$$
(\ad  X_0)^2X_1 = 0.
$$
For $n>1$, the cycle of weight 0 is
$$
p_1 \wedge q_1 + \dots + p_n\wedge q_n \eqno{(*)}
$$
and the relation expressed in terms of generators looks very ugly.  It
can be beautified as follows.  In the space of relations corresponding
to the other irreducible component, $R(\pi_2)$, the subspace of
relations of weight 0 is of dimension $n-1$.  Therefore, each summand
in $(*)$ vanishes; select the simplest one of them, say, the following
one:
$$
\{p_1, q_1\} = 0.
$$

In terms of generators this relation is:
$$
\{p_1, \{ \dots \{p_1, q_1p_2\}, q_2p_3,\} \dots, q_{n-1}p_n\},
q_n^2\}, q_{n-1}p_n\}, \dots ,\}, q_1p_2\} = 0.
$$

\ssec{3.3.  Relations in \protect $\fN_+$ for $\fh (2n)$, where
$n>1$} The space $H_1(\fn_+; \fg_1)$ is responsible for the
following relations:
$$
(\ad p_1q_2)^4(q_1^3)=0 ,  (\ad  p_2q_3)(q_1^3)=0, \dots
, (\ad  p_{n-1}q_n)(q_1^3)=0, \ad  p_n^2(q_1^3)=0.
$$
The space $H_2(\fg_+)$ is the direct sum of irreducible
$\fg_0$-modules with the following lowest weights:

\vskip 0.2 cm
\small

\begin{tabular}{|c|c|c|}
\hline
$N$ & the lowest weight & the corresponding cycle \\
\hline
1 & $-5\eps  _1-\eps  _2 $&$ q_1^3\wedge q_1^2q_2 $\\
\hline
2 &$ -3\eps  _1-3\eps  _2 $&$ 3q_1q_2^2\wedge q_1^2q_2 -q_2^3\wedge
q_1^3$\\
\hline
3 & $-2\eps  _1-2\eps  _2 $&$ \sum_i[-q_1^2q_i\wedge q_2^2p_i
+q_1^2p_i\wedge q_2^2q_i +2q_1q_2q_i\wedge q_1q_2p_i] $ \\
\hline
4 &$ -\eps  _1-\eps  _2 $&$  \sum_{i, j}[2q_1q_ip_i\wedge q_2q_ip_i
-q_1q_iq_j\wedge q_2p_ip_j +q_2q_iq_j\wedge q_1p_ip_j] $\\
\hline
5 & $0$ & $\sum_{i, j, k}[3q_iq_jp_k\wedge q_kp_ip_j
-q_iq_jq_k\wedge p_ip_jp_k] $\\
\hline
6 & $-\eps  _1$ & $q_1\sum_iq_i^3\wedge p_i^3$ \\
\hline
\end{tabular}
\normalsize
\vskip 0.2 cm

\noindent The last cocycle corresponds to the slippery relation of
degree 3.

\ssec{3.4. Relations in $\fN_+$ for $\fk (2n+1)$, where $n>1$} For
the $X_{i}^{+}$, the relations are the same as for $\fn_{+}$ of
$\fsp (2n+2)$.

The relations between $X_i^+$, where $1 \leq i\leq n$ and $Y$ are
the same as for $\fh (2n)$.

New relations, the ones involving $X_0^+$ and $Y$, are:
\vskip 0.2 cm
\small
\begin{tabular}{|c|c|c|c|}
\hline
$N$ & the lowest weight & the corresponding cycle & a simplified
relation\\
\hline
1 & $-4\eps  _1 $&$ \sum_i-q_1^2q_i\wedge q_1^2p_i
+(n+2)tq_1\wedge q_1^3 $ & \\
\hline
2 & $-3\eps  _1-\eps  _2$ &$ q_1^3 \wedge tq_2
+q_1^2q_2\wedge tq_1 $&$ [X_1^+, [Y, X_0^+]] = 0$\\
\hline
\end{tabular}
\normalsize

\ssec{3.5. Relations between $\fN_+$ and $\fN_-$ for $\fh (2n)$
and $\fk (2n+1)$, where $n>1$} These relations are as for
$\fsp(2n+2)$ unless they involve $Y$; and the new, extra, ones are
(the right hand side is just $-3q_1^2$):
$$
[Y, X_0^+] =  -\frac{3}{2^{n-1}}(\ad  X_1^+)^2 \dots (\ad
X_{n-1}^+)^2X_n^+;\quad
[Y, X_i^+] = 0\ \text{ for}\ i>0.
$$

\section*{\protect\S 4. Relations for $\fvect(n)$ and $\fsvect (m)$, where $n, m>2$}

\ssec{4.1. Relations in $\fN_-$ for $\fvect(n)$ and $\fsvect (m)$, where $n,
m>2$}
These relations   are the same as for $\fn_-$ of  $\fsl(n+1)$.

\ssec{4.2. Relations for $\fN_+$ of $\fvect (n)$, where  $n\geq
4$} The relations that  constitute $H_1(\fn_+; \fg_1)$ are of two
kinds:

1) those between $X^+_1, \dots , X^+_n$ which are the same as for $\fsl(n+1)$,
namely,
$$
\begin{matrix}(\ad  X^+_i)(X^+_j)&=&0 \; \text{ for}\; |i-j|>1,
(\ad  X^+_{i\pm 1})^2 X^+_i &=&0 .\end{matrix}
$$

2) those that involve $Y$ but not $X^+_n$:
$$
\begin{matrix}
(\ad  X^+_i)^2 Y&=&0\; \text{ for }\; i=1, n-2, n-1;\\
(\ad  X^+_i)Y &=&0 \; \text{ for}\; i\not=1, n-2, n-1.\end{matrix}
$$
The space $H_2(\fg_+)$ is  the  direct  sum  of  irreducible $\fg_0$-modules
with  the  following  lowest  weights:
\vskip 0.1 cm

\small
$$
\begin{tabular}{|c|c|c|}
\hline
$N$ & the lowest weight & the corresponding cycle  \\
\hline
1 &$ -\eps  _1-\eps  _2+4\eps  _n $&$ x_n^2\partial _1\wedge
x_n^2\partial _2 $ \\
\hline
2 &$ -2\eps  _1+\eps  _{n-1}+3\eps  _n $&$ x_n^2\partial
_1\wedge
x_nx_{n-1}\partial _1$ \\
\hline
3 &$ -\eps   _1-3\eps   _n $&$ n\sum _i(x_n^2\partial _1\wedge
x_nx_{i}\partial _i)+2\sum _i(x_n^2\partial _i\wedge
x_nx_{i}\partial _1)$  \\
\hline
4 &$ -\eps  _1-\eps  _{n-1}-2\eps  _n $&
$\sum _i(x_nx_{n-1}\partial _1\wedge
x_nx_{i}\partial _i)-\sum _i(x_n^2\partial _1\wedge
x_{n-1}x_{i}\partial _i)  $ \\
\hline
5 & $-\eps  _1-\eps  _{n-1}-2\eps  _n $&
$\sum _i(x_nx_{n-1}\partial _i\wedge
x_nx_{i}\partial _1)-\sum _i(x_n^2\partial _i\wedge
x_{n-1}x_{i}\partial _1)$ \\
\hline
6 & $-\eps  _{n-1}-\eps  _n $&
$\sum _{i, j}x_nx_{i}\partial _i\wedge
x_{n-1}x_{j}\partial _j$ \\
\hline
7 & $-\eps  _{n-1}-\eps  _n $&
$\sum _{i, j}x_nx_{i}\partial _j\wedge
x_{n-1}x_{j}\partial _i$ \\
\hline
8 & $-\eps  _{1}-\eps  _{2}+2\eps  _{n-1}+2\eps  _n $&
$2x_nx_{n-1}\partial _1\wedge
x_{n}x_{n-1}\partial _2-x_n^2\partial _1\wedge
x^2_{n-1}\partial _2-x^2_{n-1}\partial _1\wedge
x_n^2\partial _2  $ \\
\hline
\end{tabular}
$$

\normalsize
\ssec{4.3. Relations for $\fN_+$ of $\fvect (3)$}
The relations that constitute
$H_1(\fn_+; \fg_1)$ are the same as for $n>3$ except that some of the
relations distinct for large $n$ merge into one for $n=3$.

The space $H_2(\fg_+)$ is  the  direct  sum  of  irreducible $\fg_0$-modules
with  the  following  lowest  weights:
\vskip 0.2 cm

\small
$$
\begin{tabular}{|c|c|c|}
\hline
$N$ & the lowest weight & the corresponding cycle  \\
\hline
1 &$ -\eps  _1-\eps  _2+4\eps  _3 $&$ x_3^2\partial _1\wedge
x_3^2\partial _2 $ \\
\hline
2 &$ -2\eps  _1+\eps  _{2}+3\eps  _3 $&$ x_3^2\partial
_1\wedge
x_3x_{2}\partial _1$ \\
\hline
3 &$ -\eps   _1-3\eps   _3 $&$ 3\sum _ix_3^2\partial _1\wedge
x_3x_{i}\partial _i+2\sum _ix_3^2\partial _i\wedge
x_3x_{i}\partial _1$  \\
\hline
4 &$ -\eps  _1-\eps  _{2}-2\eps  _3 $&
$\sum _ix_3x_{2}\partial _1\wedge
x_3x_{i}\partial _i-\sum _ix_3^2\partial _1\wedge
x_{2}x_{i}\partial _i$ \\
\hline
5 & $-\eps  _1-\eps  _{2}-2\eps  _3 $&
$\sum _ix_3x_{2}\partial _i\wedge
x_3x_{i}\partial _1-\sum _ix_3^2\partial _i\wedge
x_{2}x_{i}\partial _1$ \\
\hline
6 & $-\eps  _{2}-\eps  _3 $&
$\sum _{i, j}x_3x_{i}\partial _i\wedge
x_{2}x_{j}\partial _j$ \\
\hline
7 & $-\eps  _{2}-\eps  _3 $&
$\sum _{i, j}x_3x_{i}\partial _j\wedge
x_{2}x_{j}\partial _i$ \\
\hline
\end{tabular}
$$

\normalsize
\ssec{4.4. The relations for $\fsvect(n)$}
These relations are those of the relations for
$\fvect(n)$ that do not involve $X^+_n =x_n\sum x_i\partial _i$.

\ssec{4.5. Relations between $\fN_+$ and $\fN_-$ for $\fvect (n)$
and $\fsvect(m)$, $n, m>2$}
These relations are as for $\fsl(n+1)$ unless they involve $Y$;
the extra relations are:
$$
[Y, X_n^+] = -2[ X_{n-1}^+, [\dots , [X_2^+, X_1^+] \dots
]](=-2x_n\partial_1),\;
[Y, X_i^+]  = 0\ \text{ for }\ i<n.
$$

\footnotesize \ssec{Acknowledgments} D.L. wishes to express his
gratitude to L.~Faddeev and A.~Kostrikin for their interest in the
preliminary versions of the manuscript (1981), to M.~Hazewinkel for
the hospitality at CWI in 1986 and 1988; to I.~Kantor, B.~Feigin and
D.~Fuchs for inspiring discussions in 1975--77, and to V.~Lebedenko,
who drew his attention to this problem.  We were financially
supported: D.L. by I.~Bendixson grant and the Swedish NFR during
1986--1990, and an NSF grant via IAS in 1989, D.L. and E.P. by
SFB--170 in June--July of 1990.  

\end{document}